\documentclass[aos,seceqn,nameyear,dvips]{arximspdf}

\doi{10.1214/09-AOS757}
\volume{38}
\issue{3}
\pubyear{2010}
\firstpage{1460}
\lastpage{1477}

\makeatletter

\newtheorem{Lemma}{Lemma}
\newproclaim{Example}{Example}
\newtheorem{Theorem}{Theorem}
\newtheorem{Proposition}{Proposition}

  \let\sv@tabnotetext\tabnotetext
  \let\sv@tabnotemark@fmt\tabnotemark@fmt
   \long\def\legend#1{{\let\tabnote@indent\leavevmode\sv@tabnotetext[]{}{#1}}}

\makeatother

\begin{document}
\begin{frontmatter}

\title{A new and flexible method for constructing designs for computer experiments}
\runtitle{Constructing designs for computer experiments}

\thankstext{t1}{Supported by grants from the Natural Sciences
and Engineering Research Council of Canada.}

\begin{aug}
\author[A]{\fnms{C. Devon} \snm{Lin}\ead[label=e1]{cdlin@mast.queensu.ca}\corref{}},
\author[B]{\fnms{Derek} \snm{Bingham}\ead[label=e2]{dbingham@stat.sfu.ca}},
\author[B]{\fnms{Randy R.} \snm{Sitter}\ead[label=e3]{sitter@stat.sfu.ca}} and
\author[B]{\fnms{Boxin} \snm{Tang}\ead[label=e4]{boxint@stat.sfu.ca}\ead[label=u1,url]{http://www.foo.com}}
\runauthor{Lin, Bingham, Sitter and Tang}
\affiliation{Queen's University, Simon Fraser University, Simon Fraser
University and~Simon~Fraser~University}
\address[A]{C. D. Lin\\
Department of Mathematics\\
\quad and Statistics\\
Jeffery Hall, University Avenue\\
Queen's University\\
Kingston, Ontario K7L 3N6\\
Canada\\
\printead{e1}}
\address[B]{D. Bingham\\
R. R. Sitter\\
B. Tang\\
Department of Statistics\\
\quad and Actuarial Science\\
8888 University Drive\\
Simon Fraser University\\
Burnaby, British Columbia V5A 1S6\\
Canada\\
\printead{e2}\\
\phantom{E-mail: }\printead*{e3}\\
\phantom{E-mail: }\printead*{e4}}
\end{aug}

% HISTORY:
\received{\smonth{8} \syear{2008}}
\revised{\smonth{9} \syear{2009}}

% ABSTRACT
%
\begin{abstract}
We develop a new method for constructing ``good'' designs for computer
experiments. The method derives its power from its basic structure that
builds large designs using small designs. We specialize the method for
the construction of orthogonal Latin hypercubes and obtain many results
along the way. In terms of run sizes, the existence problem of
orthogonal Latin hypercubes is completely solved. We also present an
explicit result showing how large orthogonal Latin hypercubes can be
constructed using small orthogonal Latin hypercubes. Another appealing
feature of our method is that it can easily be adapted to construct
other designs; we examine how to make use of the method to construct
nearly orthogonal and cascading Latin hypercubes.
\end{abstract}

% KEYWORDS
\setattribute{keyword}{AMS}{AMS 2000 subject classification.}

\end{frontmatter}

%s1 ###
\section{Introduction}\label{sec1}

Scientists are increasingly using experiments on computer
simulators to help understand physical systems. Computer experiments
differ from physical experiments in that the systems are usually
deterministic, and thus the response in computer experiments is
unchanged if a design point is replicated. The lack of random error
presents challenges to both the design and analysis of experiments
[e.g., see Sacks et al. (\citeyear{Sacksetal89})].

Similar to physical experiments, computer experiments are performed
with a variety of goals in mind. Objectives include factor screening
[Welch et al. (\citeyear{Welchetal92}), Linkletter et al.
(\citeyear{Linkletteretal06})], building an emulator of the simulator
[Sacks et al. (\citeyear{Sacksetal89})], optimization [Jones, Schonlau
and Welch (\citeyear{JSW98})] and model calibration [Kennedy and
O'Hagan (\citeyear{KOH01})]. Latin hypercube designs [McKay, Beckman
and Conover (\citeyear{MBC79})] are commonly used for computer
experiments. These designs have the feature that when projected onto
one dimension, the equally-spaced design points ensure that each of the
input variables has all portions of its range represented.

While constructing Latin hypercube designs is fairly easy, it
is more challenging to find these designs
when optimality criteria are imposed. For details of optimality
criteria, see Shewry and Wynn (\citeyear{SW87}), Morris and Mitchell (\citeyear{MM95}),
Joseph and Hung (\citeyear{JH08}) and the references therein. In this article, we
focus on the orthogonality of Latin hypercubes.
Ye (\citeyear{Y98}), Steinberg and Lin (\citeyear{SL06}) and Cioppa and
Lucas (\citeyear{CL07}) developed methods for constructing orthogonal Latin
hypercubes. These methods all have restrictions on the run size $n$.
The approach of Ye (\citeyear{Y98}) and Cioppa and Lucas (\citeyear{CL07}) gives
designs for $n = 2^k$ or $2^k+1$, and
the method of Steinberg and Lin (\citeyear{SL06}) provides
designs for $n = 2^{2^k}$
where $k \geq2$ is an integer.
Practitioners would appreciate
a methodology that can quickly produce designs with more flexible run sizes.

In this article, a new construction is proposed for finding ``good''
Latin hypercube designs for computer experiments. The method is
simple and uses small designs to construct
larger designs with desirable properties. Our
methodology is quite powerful insofar as it allows orthogonal Latin
hypercubes to be constructed for any run size $n$ where $n \neq
4k+2$. When $n=4k+2$, we prove that an orthogonal Latin hypercube
does not exist. Another important feature of our method is that it
can easily be adapted to construct
nearly orthogonal Latin hypercubes and cascading Latin hypercubes
[Handcock (\citeyear{H91})].

The article is outlined as follows. Section \ref{sec2} introduces notation,
presents a general method of construction
and discusses how to obtain Latin hypercubes based on
this general structure.
Section \ref{sec3} devotes itself to the construction
of orthogonal Latin hypercubes. Besides several general
theoretical results and many concrete examples, an existence
result is also established here.
In Section \ref{sec4}, we examine how the general method can be used
to construct nearly orthogonal Latin hypercubes.
We conclude the article with some remarks in Section \ref{sec5}.
The proofs for some theoretical results are deferred to
\hyperref[app]{Appendix}
for a smooth flow of the main ideas and results.

%s2 ###
\section{A general method of construction}\label{sec2}

Consider designs of $n$ runs with $m$ factors of $s$ levels
where $2\leq s\leq n$.
Without loss of generality,
the $s$ levels are taken to be centered at zero and equally spaced.
For odd $s$, the levels are taken as $-(s-1)/2, \ldots, -1, 0,1, \ldots,
(s-1)/2$, and for even $s$, they are $-(s-1)/2,\ldots,
-1/2,1/2,\ldots,(s-1)/2$.
The levels, except for level $0$ in the case of
odd $s$, are assumed to be equally replicated
in each design column
to ensure that linear main effects are all orthogonal to the
grand mean. Such a design is denoted by $D(n, s^m)$
and can be represented by
an $n \times m$ matrix $D = (d_{ij})$ with entries from the set of
$s$ levels as described above.
In this notation, an $m$-factor Latin hypercube
design is a $D(n,s^m)$ with $n=s$.

%s2.1 ###
\subsection{Construction method}\label{sec21}

Let $A=(a_{ij})_{n_1 \times m_1}$
be a matrix with entries $a_{ij}=\pm1$,
$B=(b_{ij})_{n_2 \times m_2}$
be a $D(n_2,s_2^{m_2})$,
$C=(c_{ij})_{n_1 \times m_1}$
be a $D(n_1,s_1^{m_1})$
and $D=(d_{ij})_{n_2 \times m_2}$ be a matrix with entries $d_{ij}=\pm1$.
Let $\gamma$ be a real number.
New designs are found using the following construction:
%
%e2.1 ###
\begin{equation}\label{eq1}
L=A\otimes B + \gamma C\otimes D,
\end{equation}
where the Kronecker product $A \otimes B$ is the $n_1n_2
\times m_1m_2$ matrix,
\[
A \otimes B = \left[
\matrix{
a_{11}B & a_{12}B & \cdots & a_{1m_1}B \cr
a_{21}B & a_{22}B & \cdots & a_{2m_1}B \cr
\vdots & \vdots & \ddots & \vdots \cr
a_{n_11}B & a_{n_12}B & \cdots & a_{n_1m_1}B}
\right]
\]
with $a_{ij}B$ itself being an $n_2 \times m_2$ matrix.
The resulting design $L$ in (\ref{eq1}) has $n=n_1n_2$ runs and
$m=m_1m_2$ factors.

The above construction has an interesting interpretation.
As an illustration, consider a simple case in which $A=(1,1)^T$ and
$C=(1/2,-1/2)^T$. Design $L$ in (\ref{eq1}) has a column,
%
%e2.2 ###
\begin{eqnarray}\label{onecolumn}
\pmatrix{b+\dfrac{\gamma}{2}d\vspace*{3pt}\cr
b-\dfrac{\gamma}{2}d},
\end{eqnarray}
where $b$ is a column of $B$ and $d$ is a column of $D$. Further let
$d = (d_1, \ldots, d_{n_2})^T$. Since $d_i = \pm1$, the
column (\ref{onecolumn}) can be viewed as simultaneously shifting
each level in $b$ to the left and the right by $\gamma/2$. If we
view $b$ as a block of level settings, then we are shifting
two identical blocks $b$, one to the left and the other to
the right. We will show in Section \ref{sec22} that with the
appropriate choices
of $A$, $B$, $C$, $D$ and~$\gamma$,
the levels in each column of $L$ in ({\ref{eq1}) are
equally spaced and unreplicated, thus resulting in a Latin hypercube.

Now consider all $m$ columns of $L$ under this simple case. Each
one-dimensional block $b$ becomes an $m$-dimensional stratum, $B$.
Suppose $D$ is a matrix of plus ones. Then the design points in $B +
\gamma D/2 $ can be obtained by shifting the entire stratum $B$ to
the right by $\gamma/2$. Similarly, the design points in
$B - \gamma D/2 $ can be obtained by shifting the entire stratum $B$
to the left by $\gamma/2$. In this case, closely
clustered points in each stratum are expected.
This feature can be utilized to construct cascading
Latin hypercubes [Lin (\citeyear{L08})].

We shall see that the orthogonality or near
orthogonality of $L$ in (\ref{eq1}) is determined by the
orthogonality or near orthogonality of $A$, $B$, $C$ and $D$, the
correlations between the columns in $A$ and those in $C$, and the
correlations between the columns in $B$ and those in $D$. As a
result, the method allows orthogonal and nearly orthogonal Latin
hypercubes to be easily constructed.

Vartak (\citeyear{V55}) appears to be the first to use the Kronecker product
systematically
to construct statistical experimental designs. In a recent work,
Bingham, Sitter and Tang (\citeyear{BST09}) introduced a method for constructing
a rich class of designs that are suitable for use in computer
experiments. Their approach corresponds to
$\gamma=0$ in the general construction given in (\ref{eq1}).
The designs in that paper have many levels
and are not Latin hypercubes in general.

%s2.2 ###
\subsection{Latin hypercubes}\label{sec22}

The following result shows how to obtain Latin hypercubes from the
construction in (\ref{eq1}).
\begin{Lemma}\label{lemma1}
Let $\gamma= n_2$. Then
design $L$ in (\ref{eq1}) is a Latin hypercube if:
\begin{longlist}[(ii)]
\item[(i)] both $B$ and $C$ are Latin hypercubes and
\item[(ii)] at least one of the following two conditions is true:
\begin{longlist}[(ii) \,(b)]
\item[(a)] $A$ and $C$ satisfy that for any $i$, if $p$ and $p'$ are
such that
$c_{pi} = - c_{p'i}$, then $a_{pi} = a_{p'i}$;
\item[(b)] $B$ and $D$ satisfy that
for any $j$, if $q$ and $q'$ are such that
$b_{qj} = - b_{q'j}$, then $d_{qj} = d_{q'j}$.
\end{longlist}
\end{longlist}
\end{Lemma}

The proof is given in the \hyperref[app]{Appendix}.
Just in terms of constructing Latin hypercubes,
Lemma \ref{lemma1} is not of much significance in itself as one can
easily obtain
a Latin hypercube simply by combining several
permutations of the set of levels. The significance of Lemma \ref{lemma1}
lies in the fact that it produces Latin hypercubes with the structure
in (\ref{eq1})
and thus provides a path to the construction of orthogonal
and cascading Latin hypercubes.

Condition (i) in Lemma \ref{lemma1} is not really a condition, and it
simply tells us to choose $B$ and $C$ to be Latin hypercubes.
In order for $L$ to be a Latin hypercube, the only mild condition
is that in (ii) of Lemma \ref{lemma1}. Two situations where condition (ii)
is obviously met are as follows:
($\alpha$) $C$ has a foldover structure
in the sense that $C =(C_0^T, -C_0^T)^T$,
and $A$ has the form $A=(A_0^T, A_0^T)^T$;
($\beta$) $A$ or $D$ is a matrix of
all plus ones. Both situations are useful.
Theorem \ref{theorem3} of Section \ref{sec33} is derived under
situation ($\alpha$).
Situation ($\beta$) can be used for constructing cascading Latin hypercubes.
We now give an example to illustrate Lemma \ref{lemma1}.
\begin{Example}\label{example1}
Consider the construction of Latin hypercubes of 32 runs with 32
factors. We choose $n_1=m_1=2$ and $n_2=m_2=16$ so that $n=n_1n_2=32$
and $m=m_1m_2=32$. To meet condition (ii) in Lemma \ref{lemma1}, let
$A$ be a matrix of all plus ones. Now let $\gamma=n_2=16$ and
$D=(d_{ij})$ be any $16 \times16$ matrix of $\pm1$. For $L$ in
(\ref{eq1}) to be a Latin hypercube, we need both $B$ and $C$ to be
Latin hypercubes. Let us use $C=[ (1/2,-1/2)^T,(-1/2,1/2)^T ]^T$ and
$B=B_0/2$ where $B_0$ is listed in Table~\ref{Table:lhd16}.
%
%t1 ###
\begin{table}
\caption{Design matrix of $B_0$ in Example \protect\ref{example1}}
\label{Table:lhd16}
\fontsize{8}{10}\selectfont{\begin{tabular*}{\tablewidth}{@{\extracolsep{\fill}}c@{}}
\\[-4pt]
\hspace*{-0.7pt}$ \left(
\begin{array}{r@{\hspace*{8pt}}r@{\hspace*{8pt}}r@{\hspace*{8pt}}r@{\hspace*{8pt}}
r@{\hspace*{8pt}}r@{\hspace*{8pt}}r@{\hspace*{8pt}}r@{\hspace*{8pt}}r@{\hspace*{8pt}}
r@{\hspace*{8pt}}r@{\hspace*{8pt}}r@{\hspace*{8pt}}r@{\hspace*{8pt}}r@{\hspace*{8pt}}r@{\hspace*{8pt}}r}
-15& 5& 9 & -3& 7& 11& -11& 7 & -9& 3& -15& 5 & 11& -11& 7& -7 \\
-13& 1& 1 & 13& -7& -11& 11& -7 & -1& -13& -13& 1 & 13& 5& 5& -3 \\
-11& 7& -7 & -11& 13& -1& -1& -13 & 9& -3& 15& -5 & -5& 11& -7& 7 \\
-9 & 3 & -15& 5& -13& 1& 1& 13 & 1& 13& 13& -1 & -13& -5& -5& 3 \\
-7 &-11& 11& -7& 11& -7& 7& 11 & 5& 15& -3& -9 & -9& 3 & 9 & 11 \\
-5 &-15& 3& 9& -11& 7& -7& -11 & 13& -1& -1& -13& -1& 9 & 11 & 15 \\
-3 &-9 & -5& -15& 1& 13& 13& -1 & -5& -15& 3& 9 & 1& 7 &-11 & -11 \\
-1 &-13& -13& 1& -1& -13& -13& 1 & -13& 1& 1& 13 & 9& -9 & -9 & -15 \\
1 & 13& 13& -1& -9& 3& -15& 5 & 11& -7& 7& 11 & -7& -7&-15 & -9 \\
3 & 9 & 5& 15& 9& -3& 15& -5 & 3& 9& 5& 15 & -15& -13&-13 & -13\\
5 & 15& -3& -9& -3& -9& -5& -15 & -11& 7& -7& -11& 15& -3& 15 & 9 \\
7 & 11& -11& 7& 3& 9& 5& 15 & -3& -9& -5& -15& 7& 15& 13 & 13 \\
9 & -3& 15& -5& -5& -15& 3& 9 & -7& -11& 11& -7 & 5 & 13 & -3 &5 \\
11 & -7& 7& 11& 5& 15& -3& -9 & -15& 5& 9& -3 & 3 & -1 & -1 &1 \\
13 & -1& -1& -13& -15& 5& 9& -3 & 7& 11& -11& 7 &-11 &-15 & 3 &-5 \\
15 & -5& -9& 3& 15& -5& -9& 3 & 15& -5& -9& 3 & -3 & 1 & 1 &-1 \\
\end{array}
\right)$
\end{tabular*}}
%
%}
%
\end{table}
According to Lemma \ref{lemma1}, design $L$ in (\ref{eq1}) is then a
$32\times 32$ Latin hypercube.
\end{Example}

%s3 ###
\section{Constructing orthogonal Latin hypercubes}\label{sec3}

We first consider in Section~\ref{sec31} the construction
of orthogonal Latin hypercubes with run sizes $n$ that are
multiples of eight. The results here are offered directly by the
construction in (\ref{eq1}). In Section \ref{sec32},
additional techniques are employed for
constructing orthogonal Latin hypercubes of other run sizes.
Results from the application of the methods in Sections \ref{sec31}
and \ref{sec32} are presented in Section \ref{sec33}.

%s3.1 ###
\subsection{Orthogonal Latin hypercubes of $n=8k$ runs}\label{sec31}

A design or matrix $X = (x_1, \ldots, x_m)$ is said to be orthogonal
if the inner product of any two columns is zero, that is, $x_i^T x_j = 0$
for all $i \neq j$. The next result provides a set of sufficient
conditions for design $L$ in (\ref{eq1}) to be orthogonal.
\begin{Lemma}\label{lemma2}
Design $L$ in (\ref{eq1}) is orthogonal if:
\begin{longlist}
\item$A$, $B$, $C$ and $D$ are all orthogonal, and
\item at least one of the two, $A^TC=0$ and $B^TD=0$, holds.
\end{longlist}
\end{Lemma}

The proof is simple, making use of the following
properties of the Kronecker product:
%
%e3.1 ###
\begin{equation}\label{eq:kron}
(A \otimes B)^T = A^T \otimes B^T
\quad\mbox{and}\quad
(A \otimes B) (C \otimes D) = (AC) \otimes(BD).
\end{equation}
Lemma \ref{lemma1} tells how to make $L$ in (\ref{eq1})
a Latin hypercube whereas Lemma \ref{lemma2} tells how to make it orthogonal.
When the two lemmas are combined, we have a way of obtaining orthogonal
Latin hypercubes.
\begin{Theorem}\label{theorem1}
Let $\gamma= n_2$. Then
design $L$ in (\ref{eq1}) is an orthogonal Latin hypercube if:

\begin{longlist}[(iii)]
\item[(i)] $A$ and $D$ are orthogonal matrices of $\pm1$;
\item[(ii)] $B$ and $C$ are orthogonal Latin hypercubes;
\item[(iii)] at least one of the two, $A^TC=0$ and $B^TD=0$, is true;
\item[(iv)] at least one of the following two conditions is true:
\begin{longlist}[(iii) \,(b)]
\item[(a)] $A$ and $C$ satisfy that for any $i$, if $p$ and $p'$ are such that
$c_{pi} = - c_{p'i}$, then $a_{pi} = a_{p'i}$;
\item[(b)]
$B$ and $D$ satisfy that
for any $j$, if $q$ and $q'$ are such that
$b_{qj} = - b_{q'j}$, then $d_{qj} = d_{q'j}$.
\end{longlist}
\end{longlist}
\end{Theorem}

The role played by $A$ and $D$ is very different from that
of $B$ and $C$ in Theorem~\ref{theorem1}. To help understand Theorem
\ref{theorem1},
one may think that $B$ and $C$ are the building material
while $A$ and $D$ provide a blueprint for the construction.
Small orthogonal Latin hypercubes $B$ and $C$ are
used to construct a large orthogonal Latin hypercube $L$
in Theorem \ref{theorem1}. Exactly how the construction is accomplished
is guided
by $A$ and $D$ which are orthogonal matrices of $\pm1$.
In addition to the right blueprint and building material, a considerable
amount of care is necessary for the final structure to be right.
This is achieved via $\gamma= n_2$ and conditions (iii) and (iv)
in Theorem \ref{theorem1}.

Note that $A$ and $D$ may or may not be square matrices, and
the orthogonality of $A$ and $D$ is imposed on their columns.
In some mathematics literature, such matrices are called Hadamard
submatrices. For convenience, we simply call $A$ or $D$ an
orthogonal matrix when its columns are orthogonal.
Hadamard matrices and orthogonal arrays with levels $\pm1$
are all such orthogonal matrices in our terminology.
A Hadamard matrix is a square orthogonal matrix of $\pm1$.
An orthogonal array with two levels $\pm1$ requires that
each of the four combinations
$(-1, -1)$,
$(-1, +1)$,
$(+1, -1)$ and
$(+1, +1)$
occurs the same number of times in every two columns.
For some comprehensive discussion on these
and other topics in the theory of factorial designs,
we refer to Dey and Mukerjee (\citeyear{DM99}), Hedayat, Sloane
and Stufken (\citeyear{HSS99}) and Mukerjee and Wu (\citeyear{MW06}).

Because of the orthogonality of $A$ and $D$, we must have
that $n_1$ and $n_2$ are equal to two or multiples of four.
The case where $n_1 = n_2 = 2$ is trivial.
Consequently, Theorem \ref{theorem1} can be used to
construct orthogonal Latin hypercubes of $n=8k$ runs, thereby
providing designs that are unavailable in Ye (\citeyear{Y98}) and Steinberg and Lin
(\citeyear{SL06}).
When $n = n_1 n_2$ is a multiple of 16,
Theorem \ref{theorem1} becomes more powerful. This point will be highlighted
in Section \ref{sec33}. We now revisit Example \ref{example1}
for an illustration of Theorem \ref{theorem1}.
\begin{Example}\label{example2}
In Example \ref{example1}, the first 12 columns of $B$ form
a $16$-run orthogonal Latin hypercube constructed by Steinberg and Lin
(\citeyear{SL06}). If $D$ is chosen to be a Hadamard matrix of order 16
in Example \ref{example1}, Theorem \ref{theorem1} tells us the first 12
columns of
$L$ in Example \ref{example1} constitute a $32 \times12$
orthogonal Latin hypercube which has one more orthogonal factor than
the $32 \times11$ orthogonal Latin hypercube obtained
by Cioppa and Lucas (\citeyear{CL07}).

When $n_1 = n_2$, a stronger result than Theorem \ref{theorem1} can be
established,
again using the properties of the Kronecker product given in (\ref{eq:kron}).
\end{Example}
\begin{Proposition}\label{proposition1}
If $n_1=n_2=n_0$ and $A, B, C, D$ and $\gamma$
are chosen according to Theorem \ref{theorem1},
then design ($L, U$) is an orthogonal Latin hypercube with $2m_1m_2$
factors where $L$ is as in Theorem \ref{theorem1} and
$U=-n_0 A \otimes B + C \otimes D$.
\end{Proposition}

We now discuss how to choose $A, B, C, D$ and $\gamma$
to construct orthogonal Latin hypercubes.
According Theorem \ref{theorem1}, we have that $\gamma= n_2$.
Matrices $A$ and $D$ need to be orthogonal with entries of $\pm1$.
As discussed earlier, two level orthogonal arrays and
Hadamard matrices are all such orthogonal matrices.
Theorem \ref{theorem1} requires that designs $B$ and $C$ be orthogonal
Latin hypercubes.
All known orthogonal Latin hypercubes from the existing literature can
be used here.
Later in this paper (see Table \ref{Table:m}), we obtain a collection of
small orthogonal Latin hypercubes through a computer search for this purpose.
So far, all are straightforward. The nontrivial aspect from applying
Theorem \ref{theorem1}
is to satisfy conditions (iii) and (iv) which require that
$A$ and $C$ (or $B$ and $D$) jointly have certain properties.
In this paper, we satisfy these two conditions by choosing
$A$ of form $A = (A_0^T, A_0^T)^T$ and $C$ of form
$C = (C_0^T, -C_0^T)^T$ where $A_0$ and $C_0$ are such that
all the columns in the matrix,
%
%e3.2 ###
\begin{equation}\label{eq:ac}
(A, C) = \left[\matrix{
A_0 & C_0 \cr
A_0 & -C_0}
\right],
\end{equation}
are mutually orthogonal.
In Section \ref{sec33} we provide a method of finding such orthogonal
matrices with the structure in (\ref{eq:ac}) when proving
Theorem \ref{theorem3}. Comments similar to those in this paragraph
can also be made regarding the application of Proposition \ref
{proposition2} in Section \ref{sec32}.

%s3.2 ###
\subsection{Orthogonal Latin hypercubes with other run sizes}\label{sec32}

Consider an orthogonal Latin hypercube of $n$ runs
with $m \geq2$ factors.
Trivially, run size $n$ cannot be two or three.
So we must have $n \geq4$. The next result
provides a complete characterization of the existence
of an orthogonal Latin hypercube in terms of run size~$n$.
\begin{Theorem}\label{theorem2}
There exists an orthogonal Latin hypercube of $n \geq4$
runs with more than one factor
if and only if $n \neq4k+2$ for any integer $k$.
\end{Theorem}

The \hyperref[app]{Appendix} contains a proof for Theorem \ref
{theorem2}. Equivalently, Theorem \ref{theorem2}
says that the run size of an orthogonal Latin hypercube has to be
odd or a multiple of 4. Theorem \ref{theorem1} provides a method for
constructing orthogonal Latin hypercubes of $n=8k$ runs.
The present section examines how
to construct orthogonal Latin hypercubes of other run sizes.

The basic idea of our method is quite simple. To obtain
an orthogonal Latin hypercube, we stack up two
orthogonal designs with mutually exclusive and exhaustive
sets of levels. To make it precise,
we use $\mathcal{S}$ to denote the set of $n$ levels of a Latin
hypercube of
$n$ runs.
Let $\mathcal{S}=\mathcal{S}_a \cup\mathcal{S}_b$ where
$\mathcal{S}_a \cap\mathcal{S}_b = \phi$, and let $n_a$ and $n_b$ be the
numbers of levels in $\mathcal{S}_a$ and $\mathcal{S}_b$,
respectively. Suppose that there exist an $n_a \times m$ orthogonal
design $D_a$ with levels in $\mathcal{S}_a$ and an $n_b\times m$
orthogonal design $D_b$ with levels in $\mathcal{S}_b$, where
for both $D_a$ and $D_b$, each level appears precisely
once within each column. Then
%
%e3.3 ###
\begin{equation}\label{eq2} L = \pmatrix{
D_a\cr
D_b}
\end{equation}
is an $n \times m$ orthogonal Latin hypercube with $n = n_a + n_b$.
Note that $D_a$ and $D_b$ themselves are not necessarily Latin hypercubes.

We consider two special choices for $\mathcal{S}_a$ and
$\mathcal{S}_b$. For easy reference later in the paper,
we call them two stacking methods.
Our first stacking method chooses $n_a$ and $n_b$ such that $|n_a-n_b|=1$
with the corresponding
$\mathcal{S}_a=\{-(n_a-1), -(n_a-3), \ldots, n_a-3,n_a-1\}$ and
$\mathcal{S}_b=\{-(n_b-1), -(n_b-3), \ldots,\break n_b-3,n_b-1\}$. This
implies that both $D_a/2$ and $D_b/2$ in (\ref{eq2}) are orthogonal
Latin hypercubes.
We may assume that $n_a$ is odd and $n_b$ is even in the above.
By Theorem \ref{theorem2}, we know that $n_b$ has form $4k$.
It follows that $n_a$ has form $4k-1$ or $4k+1$. Thus
the first stacking method allows
orthogonal Latin hypercubes of run sizes $8k-1$ and $8k+1$ to be
constructed.

The second stacking method is more generally applicable
and it chooses $\mathcal{S}_a=\{-(n_a-1)/2, -(n_a-3)/2,
\ldots, (n_a-3)/2,(n_a-1)/2\}$ and
%
%e3.4 ###
\begin{equation}\label{eq:sb}
\mathcal{S}_b=\{-(n-1)/2,
\ldots, -(n_a+1)/2, (n_a+1)/2, \ldots, (n-1)/2\},
\end{equation}
where $n=n_a+n_b$. For this
choice, $D_a$ is an orthogonal Latin hypercube while $D_b$ is not.
We examine how to construct an orthogonal design $D_b$
with level set
$\mathcal{S}_b$ given in (\ref{eq:sb}).
Now consider the matrices in Table \ref{Table:od}.
Each of the four matrices in Table \ref{Table:od} has the following
properties:
(i) it has real entries $\pm x_1, \ldots, \pm x_{n/2}$;
(ii)~both $x_i$ and $-x_i$ occur exactly once in each column;
(iii) every two columns are orthogonal.
We note that the matrices in Table \ref{Table:od} are related to but
different from orthogonal designs in the combinatorics literature
[Geramita and Seberry (\citeyear{GS79})].

%
%t2 ###
\begin{table}
\tabcolsep=0pt
\caption{Four useful matrices}\label{Table:od}
\fontsize{9.5}{12}\selectfont
{\begin{tabular*}{\tablewidth}{@{\extracolsep{\fill}}lcccccccccccccc@{}}
\hline
\multicolumn{15}{c}{$\bolds n$}\\
\hline
\textbf{2} & \multicolumn{2}{c}{\textbf{4}} &\multicolumn{4}{c}{\textbf{8}}&
\multicolumn{8}{c@{}}{\textbf{16}} \\
\hline
$x_1$  &  $x_1$    & $x_2$   &  $x_1$    & $-x_2$  & $x_4$   & $x_3$    & $x_1$    & $-x_2$    & $-x_4$   & $-x_3$   & $-x_8$   & $x_7$  & $x_5$  & $x_6$  \\
$-x_1$ &  $x_2$    & $-x_1$  &  $x_2$    & $x_1$   & $x_3$   & $-x_4$   & $x_2$    & $x_1$     & $-x_3$   & $x_4$    & $-x_7$   & $-x_8$ & $-x_6$ & $x_5$  \\
       &  $-x_1$   & $-x_2$  &  $x_3$    & $-x_4$  & $-x_2$  & $-x_1$   & $x_3$    & $-x_4$    & $x_2$    & $x_1$    & $-x_6$   & $-x_5$ & $x_7$  & $-x_8$ \\
       &  $-x_2$   & $x_1$   &  $x_4$    & $x_3$   & $-x_1$  & $x_2$    & $x_4$    & $x_3$     & $x_1$    & $-x_2$   & $-x_5$   & $x_6$  & $-x_8$ & $-x_7$ \\
       &           &         &  $-x_1$   & $x_2$   & $-x_4$  & $-x_3$   & $x_5$    & $-x_6$    & $-x_8$   & $x_7$    & $x_4$    & $x_3$  & $-x_1$ & $-x_2$ \\
       &           &         &  $-x_2$   & $-x_1$  & $-x_3$  & $x_4$    & $x_6$    & $x_5$     & $-x_7$   & $-x_8$   & $x_3$    & $-x_4$ & $x_2$  & $-x_1$ \\
       &           &         &  $-x_3$   & $x_4$   & $x_2$   & $x_1$    & $x_7$    & $-x_8$    & $x_6$    & $-x_5$   & $x_2$    & $-x_1$ & $-x_3$ & $x_4$  \\
       &           &         &  $-x_4$   & $-x_3$  & $x_1$   & $-x_2$   & $x_8$    & $x_7$     & $x_5$    & $x_6$    & $x_1$    & $x_2$  & $x_4$  & $x_3$  \\
       &           &         &           &         &         &          & $-x_1$   & $x_2$     & $x_4$    & $x_3$    & $x_8$    & $-x_7$ & $-x_5$ & $-x_6$ \\
       &           &         &           &         &         &          & $-x_2$   & $-x_1$    & $x_3$    & $-x_4$   & $x_7$    & $x_8$  & $x_6$  & $-x_5$ \\
       &           &         &           &         &         &          & $-x_3$   & $x_4$     & $-x_2$   & $-x_1$   & $x_6$    & $x_5$  & $-x_7$ & $x_8$  \\
       &           &         &           &         &         &          & $-x_4$   & $-x_3$    & $-x_1$   & $x_2$    & $x_5$    & $-x_6$ & $x_8$  & $x_7$  \\
       &           &         &           &         &         &          & $-x_5$   & $x_6$     & $x_8$    & $-x_7$   & $-x_4$   & $-x_3$ & $x_1$  & $x_2$  \\
       &           &         &           &         &         &          & $-x_6$   & $-x_5$    & $x_7$    & $x_8$    & $-x_3$   & $x_4$  & $-x_2$ & $x_1$  \\
       &           &         &           &         &         &          & $-x_7$   & $x_8$     & $-x_6$   & $x_5$    & $-x_2$   & $x_1$  & $x_3$  & $-x_4$ \\
       &           &         &           &         &         &          & $-x_8$   & $-x_7$    & $-x_5$   & $-x_6$   & $-x_1$   & $-x_2$ & $-x_4$ & $-x_3$ \\

\hline
\end{tabular*}}
\end{table}

%
%t3 ###
\begin{table}[b]
\caption{The maximum number $m$ of columns in $\operatorname{OLH}(n,m)$
by the algorithm for $4 \leq n \leq21$}
\label{Table:m}
\begin{tabular*}{\tablewidth}{@{\extracolsep{\fill}}lcccccccccccccc@{}}
\hline
$n$ & 4 & 5& 7 & 8 & 9 & 11 & 12 & 13 & 15 & 16 & 17 & 19 & 20 & 21\\
[4pt]
$m$ & 2 & 2 & 3 & 4 & 5 & 7 & 6 & 6 & 6 & 12 & 6 & 6 & 6 & 6 \\
\hline
\end{tabular*}
\end{table}

The matrices in Table \ref{Table:od} can be used to
construct orthogonal Latin hypercubes of $n$ runs by setting
$x_i = (2i -1)/2$ for $i =1, \ldots, n/2$.
They also provide a direct construction of orthogonal designs $D_b$
with level set $\mathcal{S}_b$ in (\ref{eq:sb})
by choosing $x_i = (n_a + 2i -1)/2$ for $i=1, \ldots, n_b/2$.
Most importantly, they are useful in the
following result that allows us to construct
$D_b$ with level set
$\mathcal{S}_b$ in (\ref{eq:sb})
for more general $n_b$.
\begin{Proposition}\label{proposition2}
Let $\gamma= 1$. Then
design $L$ in (\ref{eq1}) is an orthogonal design with level set
$\{-(n_a+n-1)/2,
\ldots, -(n_a+1)/2, (n_a+1)/2, \ldots, (n_a+n-1)/2\}$ if:

\begin{longlist}
\item[(i)] $A$ and $D$ are orthogonal matrices of $\pm1$;
\item[(ii)] $B$ is an orthogonal Latin hypercube, and $C$ is an orthogonal
design with level set
$\pm(n_a+ n_2 )/2,
\pm(n_a+ 3 n_2 )/2,
\ldots,
\pm(n_a+ (n_1- 1) n_2 )/2$;
\item[(iii)] at least one of the two, $A^TC=0$ and $B^TD=0$, is true;
\item[(iv)] at least one of the following two conditions is true:
\begin{longlist}[(b) \,(iii)]
\item[(a)] $A$ and $C$ satisfy that for any $i$, if $p$ and $p'$ are such that
$c_{pi} = - c_{p'i}$, then $a_{pi} = a_{p'i}$;
\item[(b)]
$B$ and $D$ satisfy that
for any $j$, if $q$ and $q'$ are such that
$b_{qj} = - b_{q'j}$, then $d_{qj} = d_{q'j}$.
\end{longlist}
\end{longlist}
\end{Proposition}

Orthogonality of design $L$ follows from Lemma \ref{lemma2}. That $L$
has a desired
set of levels can easily be established which follows a similar path
to that for Lemma \ref{lemma1}.
Comparing Proposition \ref{proposition2} with Theorem \ref{theorem1},
we see that the only changes are
those made to $\gamma$ and $C$. Mathematically, Theorem \ref{theorem1} is
a special case of Proposition \ref{proposition2} as one can obtain the former
from the latter by setting $n_a = 0$. We present them separately
because they carry different messages and serve different purposes
in this paper.

Design $C$ required in Proposition \ref{proposition2}
can easily be obtained from the matrices in Table \ref{Table:od}.
By letting $n = n_b$ in Proposition \ref{proposition2}, design $L$ in
Proposition \ref{proposition2}
can then used as our $D_b$ as it has desired level set
$\mathcal{S}_b$ in (\ref{eq:sb}).
The run size $n_b$ of such $D_b$ has form $n_b = 8k$.
Since there is no restriction in the run size $n_a$ of $D_a$, other
than that $D_a$ is an orthogonal Latin hypercube, this second
stacking method allows orthogonal Latin hypercubes of any run size
$n \neq4k+2$ to be constructed.
\begin{Example}\label{example3}
In Example \ref{example2}, if we choose $\gamma=1$ and
let $C = (-17/2, 17/2)^{T}$, Proposition \ref{proposition2}
gives an orthogonal design $D_b$ of $n_b = 32$ runs
for 12 factors, where each column of $D_b$ is a permutation of
$-16, -15, \ldots, -1$, $1, \ldots, 15, 16$.
Now let $n_a=1$ and $D_a$ be a row of zeros. Then stacking
up $D_a$ and $D_b$ gives a $33 \times12$ orthogonal Latin
hypercube.
\end{Example}

%s3.3 ###
\subsection{Some results}\label{sec33}

The methods in Sections \ref{sec31} and \ref{sec32} both build large
orthogonal Latin hypercubes from small ones. To apply the methods, we
need to find orthogonal Latin hypercubes with small runs. Various
efficient algorithms can be helpful in this regard. Lin
(\citeyear{L08}) reported an algorithm adapted from that of Xu
(\citeyear{X02}). The key idea of the algorithm is to add columns
sequentially to an existing design. To add a column, two operations,
pairwise switch and exchange, are used. A pairwise switch switches a
pair of distinct levels in a column. For a candidate column, the
algorithm searches for all possible pairwise switches and makes the
pairwise switch that achieves the best improvement. This search and
pairwise switch procedure is repeated until an orthogonal Latin
hypercube is found. An exchange replaces the candidate column by a
randomly generated column. The exchange step is repeated at most $T_1$
(user-specified) times if no orthogonal Latin hypercube is obtained.
Since the procedure relies on the initial random columns, the entire
procedure is repeated $T_2$ times. Apart from the sequential idea, the
efficiency of the algorithm benefits from its fast updates of
orthogonality. An update is needed when a pairwise switch is applied.
The maximum number $m$ of the columns in orthogonal Latin hypercubes of
$n$ runs found by the algorithm is given in Table~\ref{Table:m} for $4
\leq n \leq21$ except for $n=16$, in which case, our algorithm finds
$m=6$. The entry $m=12$ for $n=16$ in Table \ref{Table:m} is due to
Steinberg and Lin (\citeyear{SL06}). The detailed design matrices for
the orthogonal Latin hypercubes in Table \ref{Table:m} are presented in
Lin (\citeyear{L08}) and also available from the authors.

For a concise presentation of the results in this section,
we use $\operatorname{OLH}(n,m)$ to denote an orthogonal Latin hypercube
of $n$ runs for $m$ factors.
We now present a general result from the application of
Theorem \ref{theorem1} in Section \ref{sec31} and the second stacking
method in Section \ref{sec32}.
\begin{Theorem}\label{theorem3}
Suppose that an $\operatorname{OLH}(n,m)$ is available
where $n$ is a multiple of 4 such that a Hadamard matrix
of order $n$ exists. Then we have that:

\begin{longlist}
\item the following orthogonal Latin hypercubes,
an $\operatorname{OLH}(2n,m)$,
an $\operatorname{OLH}(4n$,\break $2m)$,
an $\operatorname{OLH}(8n,4m)$ and
an $\operatorname{OLH}(16n,8m)$,
can all be constructed;
\item
all the following orthogonal Latin hypercubes,
an $\operatorname{OLH}(2n+1,m)$,
an $\operatorname{OLH}(4n+1,2m)$,
an $\operatorname{OLH}(8n+1,4m)$ and
an $\operatorname{OLH}(16n+1,8m)$
can also be constructed.
\end{longlist}
\end{Theorem}

We give a proof for Theorem \ref{theorem3}. The proof in fact provides
a detailed procedure for the actual construction of
these orthogonal Latin hypercubes.
Part (i) of Theorem \ref{theorem3} results from an application
of Theorem \ref{theorem1} in Section \ref{sec31}.
In the general construction (\ref{eq1}), we choose $B$ to be the given
$\operatorname{OLH}(n,m)$.
Matrix $D$ is obtained by taking $m$ columns from a Hadamard
matrix of order $n$. Design $C$ is chosen to be an orthogonal
Latin hypercube derived from a matrix in Table \ref{Table:od}. Note that
each of the four matrices in Table \ref{Table:od} has a fold-over structure
in that it can be written as $(X^T, -X^T)^T$. Now let
$A = (S^T, S^T)^T$ where $S$ is
obtained from $X$ by setting $x_i = 1$ for all $i$.
With the above choices for $A$, $B$, $C$ and $D$,
conditions (i), (ii), (iii) and (iv) in Theorem \ref{theorem1}
are all satisfied. This proves part (i) of Theorem \ref{theorem3}.
The proof for part (ii) of Theorem \ref{theorem3} is similar, involving
the second stacking method with $n_a = 1$ and an application
of Proposition \ref{proposition2}.

Theorem \ref{theorem3} is a very powerful result. By repeated application
of Theorem \ref{theorem3}, one can obtain many infinite series of orthogonal
Latin hypercubes. For example,
starting with an $\operatorname{OLH}(12,6)$ from Table \ref{Table:m},
we can obtain
an $\operatorname{OLH}(192,48)$ which can be used in turn to construct
an $\operatorname{OLH}(768,96)$ and so on. For another example,
an $\operatorname{OLH}(256,248)$ in Steinberg and Lin (\citeyear{SL06})
can be used to construct
an $\operatorname{OLH}(1024,496)$, an $\operatorname{OLH}(4096,1984)$
and so on.

One important problem in the study of orthogonal
Latin hypercubes is to determine the maximum
number $m^*$ of factors for an $\operatorname{OLH}(n,m^*)$
to exist. Theorem \ref{theorem2} says that $m^*=1$ if $n$ is 3 or
has form $n=4k+2$ and that $m^* \geq2$ otherwise.
This result is now strengthened below.
\begin{Proposition}\label{proposition3}
The maximum number $m^*$ of factors for an orthogonal
Latin hypercube of $n = 16k + j$ runs
has a lower bound given below:
\begin{longlist}
\item  $m^* \geq6$ for all $n=16k+j$ where $k \geq1$ and
$j \neq2, 6, 10, 14$;
\item $m^* \geq7$ for $n=16k + 11$ where $k \geq0$;
\item $m^* \geq12$ for $n=16k, 16k + 1$ where $k \geq2$;
\item $m^* \geq24$ for $n=32k, 32k + 1$ where $k \geq2$;
\item $m^* \geq48$ for $n=64k, 64k + 1$ where $k \geq2$.
\end{longlist}
\end{Proposition}

Part (i) of Proposition \ref{proposition3} is obtained as follows.
By our second stacking method with the use
of the $16 \times8$ matrix in Table \ref{Table:od}, we can construct
an $\operatorname{OLH}(n+16, m)$ where $m \leq8$ if an $\operatorname
{OLH}(n,m)$ is
available. Part (i) of Proposition \ref{proposition3} will be true if
we can claim that an $\operatorname{OLH}(n,6)$ exists for
all $ 17 \leq n \leq32 $ except for $n=18,22, 26, 30$.
We already know that the claim is true for $n=17, 19, 20, 21$
from Table \ref{Table:m} and for $n=32$ from Example \ref{example2}.
Note that an $\operatorname{OLH}(11,6)$ can be obtained by choosing any
six columns
from the $\operatorname{OLH}(11,7)$ in Table \ref{Table:m}.
For $n=23$, we use the first stacking method
by choosing $n_a = 11$ and $n_b = 12$ and using an $\operatorname
{OLH}(11,6)$ and
the $\operatorname{OLH}(12,6)$ in Table \ref{Table:m}.
The case
$n=24$ follows from applying part (i) of Theorem \ref{theorem3}
to the $\operatorname{OLH}(12,6)$ in Table \ref{Table:m}.
For $n=25$, an $\operatorname{OLH}(25,6)$ can be constructed using
the first stacking method with $n_a = 13$ and $n_b = 12$.
For $n=27$, we apply the second stacking method
by choosing $n_a = 11$ and $n_b = 16$.
The second stacking method also allows the construction
of an $\operatorname{OLH}(28,6)$, an $\operatorname{OLH}(29,6)$ and an
$\operatorname{OLH}(31,6)$.
We choose $n_a = 12$ and $n_b =16$ for $n=28$,
$n_a = 13$ and $n_b = 16$ for $n=29$, and
$n_a = 15$ and $n_b =16$ for $n=31$.
Part (ii) follows from the existence of an $\operatorname{OLH}(11,7)$
in Table \ref{Table:m}.
Parts (iii), (iv) and (v)
follows from an application of Theorem \ref{theorem3}.

The following remarks are in order regarding Proposition \ref{proposition3}.
If we wish, we can obtain sharper lower bounds
on $m^*$ for certain values of $n$ by applying Theorem \ref{theorem3}.
For example, using the $\operatorname{OLH}(12,6)$ in Table \ref{Table:m},
we can establish that $m^* \geq6 \times8^k$ for $n=12 \times16^k$.
We will not dwell further on this issue but are satisfied
with the general lower bound in Proposition \ref{proposition3}.
The lower bound in Proposition \ref{proposition3} is derived
from the small orthogonal Latin hypercubes found by our
algorithm. Therefore, improved bounds will be naturally
available in the future
if better results are obtained from computer search.

%
%t4 ###
\begin{table}[b]
\tablewidth=250pt
\caption{Orthogonal Latin hypercubes of
$n=16k$ runs where $k \geq2$}
\label{Table:comp}
\begin{tabular*}{\tablewidth}{@{\extracolsep{\fill}}lcrrc@{}}
\hline
$\bolds n$ & \multicolumn{1}{c}{$\bolds m$} & \multicolumn{1}{c}{\textbf{Ye}}
& \multicolumn{1}{c}{\textbf{SL}} & \multicolumn{1}{c@{}}{\textbf{CL}}\\
\hline
\phantom{0}32 & 12 & 8 & 0 & 11\\
\phantom{0}48 & 12 & 0 & 0 & \phantom{0}0 \\
\phantom{0}64 & \phantom{0}32$^*$ & 10& 0 & 16\\
\phantom{0}80 & 12 & 0 & 0 & \phantom{0}0 \\
\phantom{0}96 & 24 & 0 & 0 & \phantom{0}0 \\
112 & 12 & 0 & 0 & \phantom{0}0\\
128 & 48 & 12& 0 & 22\\
144 & \phantom{0}24$^*$ & 0 & 0 & \phantom{0}0\\
160 & 24 & 0 & 0 & \phantom{0}0 \\
176 & 12 & 0 & 0 & \phantom{0}0\\
192 & 48 & 0 & 0 & \phantom{0}0\\
208 & 12 & 0 & 0 & \phantom{0}0 \\
224 & 24 & 0 & 0 & \phantom{0}0 \\
240 & 12 & 0 & 0 & \phantom{0}0\\
256 & 192$^*$ & 14 & 248 & 29\\
\hline
\end{tabular*}
\legend{Note: Ye: the number of orthogonal columns by Ye
(\citeyear{Y98}); SL: the number of orthogonal columns by Steinberg and
Lin (\citeyear{SL06}); CL: the number of orthogonal columns by Cioppa
and Lucas (\citeyear{CL07}).}
\end{table}

Lin (\citeyear{L08}) in her thesis provides a comprehensive table
of orthogonal Latin hypercubes for all $n \leq256$.
Here we present the results in Table \ref{Table:comp} for the case where
$n$ is a multiple of 16.
The first column is the run size and the
second column is the number of factors
obtained by our methods. Those entries marked with
an $*$ are given by Proposition \ref{proposition1}.
The remaining columns
of Table \ref{Table:comp} give the number of factors obtained by
the methods of Ye (\citeyear{Y98}), Steinberg and Lin (\citeyear{SL06})
and Cioppa and Lucas (\citeyear{CL07}).
Table \ref{Table:comp} clearly shows that
our methods can provide orthogonal Latin hypercubes when
other methods cannot be applied. When
other methods are applicable, our
methods give many more factors than these existing methods
with the only exception given
by $n=256$, for which case Steinberg and Lin (\citeyear{SL06}) found
an $\operatorname{OLH}(256,248)$.

%s4 ###
\section{Nearly orthogonal Latin hypercubes}\label{sec4}

The general construction in (\ref{eq1}) is very versatile and can also be
used to construct nearly orthogonal and cascading Latin hypercubes.
Due to space limitation, we omit the discussion on cascading Latin hypercubes
and refer the reader to Lin (\citeyear{L08}). In what follows,
we provide a brief discussion on nearly orthogonal Latin hypercubes;
interested readers can find more details in Lin's thesis (\citeyear{L08}).

To assess near orthogonality, we
adopt two measures defined in Bingham, Sitter and Tang (\citeyear{BST09}). For a
design $D=(d_1, \ldots, d_m)$, where $d_j$ is the $j$th column of~$D$, define $\rho_{ij}(D)$ to be
$d_i^T d_j/[d_i^T d_i d_j^Td_j]^{1/2}$.
If the mean of the level settings in $d_j$
for all $j=1,\ldots,m$ is zero, then $\rho_{ij}(D)$ is simply the
correlation coefficient between columns $d_i$ and $d_j$.
Near orthogonality can be
measured by the maximum correlation
$\rho_M(D)=\max_{i,j}|\rho_{ij}(D)|$ and the average
squared correlation $\rho^2(D)=\sum_{i<j}\rho^2_{ij}(D)/[(m(m-1)/2]$.
Smaller values of $\rho_M(D)$ and $\rho^2(D)$ imply near
orthogonality. Obviously, if $\rho_M(D)=0$ or $\rho^2(D)=0$, then
an orthogonal Latin hypercube is obtained. The following result
shows how the method in (\ref{eq1}) can be used to construct
nearly orthogonal Latin hypercubes.
\begin{Proposition}\label{proposition4}
Suppose that $A$, $B$, $C$, $D$ and $\gamma$ in (\ref{eq1}) are
chosen according to Lemma \ref{lemma1} so
that design $L$ in (\ref{eq1}) is a Latin hypercube.
In addition, we assume that $A$ and $D$ are orthogonal and
that at least one of the two, $A^TC = 0$ and $B^TD = 0$,
holds true. We then have that:
\begin{longlist}
\item $\rho^2(L)=w_1\rho^2(B)+w_2\rho^2(C)$;
\item $\rho_M(L)=\operatorname{Max}\{w_3\rho_M(B),
w_4\rho_M(C)\}$,
\end{longlist}
where $w_1$, $w_2$, $w_3$ and $w_4$ are given by
$w_1=(m_2-1)(n_2^2-1)^2/[(m_1m_2-1)(n^2-1)^2]$,
$w_2=n_2^4(m_1-1)(n_1^2-1)^2/[(m_1m_2-1)(n^2-1)^2]$,
$w_3=(n_2^2-1)/(n^2-1)$ and
$w_4=n_2^2(n_1^2-1)/(n^2-1)$.
\end{Proposition}

The proof for Proposition \ref{proposition4} is in the \hyperref[app]{Appendix}.
Proposition \ref{proposition4} says that if $B$ and $C$ are nearly orthogonal,
the resulting Latin hypercube $L$ is also nearly orthogonal.
An example, illustrating
the use of this result, is considered below.

%
%t5 ###
\begin{table}
\caption{Design matrix of $B_0$ in Example \protect\ref{example4}}
\label{Table:nolhd16}
\fontsize{8.5}{10}\selectfont{\begin{tabular*}{\tablewidth}{@{\extracolsep{\fill}}c@{}}
\\[-4pt]
\hspace*{4.3pt}$ \left(
\begin{array}{r@{\hspace*{8pt}}r@{\hspace*{8pt}}r@{\hspace*{8pt}}r@{\hspace*{8pt}}
r@{\hspace*{8pt}}r@{\hspace*{8pt}}r@{\hspace*{8pt}}r@{\hspace*{8pt}}r@{\hspace*{8pt}}
r@{\hspace*{8pt}}r@{\hspace*{8pt}}r@{\hspace*{8pt}}r@{\hspace*{8pt}}r@{\hspace*{8pt}}r@{}}
-15& 15& -13& 13& -5& -13& 5& 3& -1& 5 & -7 & 5& -9 & -9&
5\\
-13& -15& -3& 3& 7& 3& 15& -11& 13& -5 & 7 &-13& -7 & -3& -3 \\
-11& -9& -5& -11& -15& 13& -5& 11& -9& 9 & 9 & 3& -5 & -1& -11 \\
-9& -1& 9& -15& -11& 1& -1& -13& 5& -1 &-15 & 7& 1 & 3& 15 \\
-7& 1& -7& 7& 15& 15& -13& 9& -5& -13 & -3 & -1& -1 & 7& 13 \\
-5& 13& 11& -5& 9& -7& -3& -9& -13& 11 & 13 & -9& -3 & 13& 1 \\
-3& -5& 13& 15& -9& -9& -11& 1& 7& -9 & 15 & 11& 9 & 1& -1 \\
-1& -11& 3& -7& 11& -15& 13& 15& -7& -3 & -9 & 9& 7 & 9& -5 \\
1& 3& -9& -3& -1& -5& -15& -1& 11& 3 &-11 &-15& 15 & 5& -15 \\
3& -3& 15& 11& 3& 9& 1& -7& -15& 1 &-13 & -3& 3 &-15& -9 \\
5& 9& 7& -1& 5& 11& 9& 13& 15& 15 & 5 & 1& 11 & -7& 9 \\
7& 7& -1& -13& 13& -1& -7& -5& 9& -7 & 3 & 15& -13 &-11& -13 \\
9& 5& -11& -9& -7& -3& 7& -3& -11& -15 & 11 & -7& 13 &-13& 7 \\
11& 11& 5& 5& -13& 7& 11& 5& 3& -11 & -5 & -5& -11 & 15& -7 \\
13& -7& -15& 9& 1& 5& 3& -15& -3& 13 & 1 & 13& 5 & 11& 3 \\
15& -13& 1& 1& -3& -11& -9& 7& 1& 7 & -1 &-11& -15 & -5& 11
\end{array}
\right)$
\end{tabular*}}
\end{table}

\begin{Example}\label{example4}
Let $A=(1,1)^{T}$, $C=(1/2,-1/2)^T$, and $\gamma=16$. Choose a $16
\times15$ nearly orthogonal Latin hypercube $B=B_0/2$ where $B_0$ is
displayed in Table \ref{Table:nolhd16},
and $B$ has $\rho^2(B)=0.0003$ and $\rho_M(B)=0.0765$. Taking any
15 columns of a Hadamard matrix of order 16 to be $D$ and then applying
(\ref{eq1}), we obtain a Latin hypercube $L$ of 32 runs
and 15 factors. As $\rho^2(C)=\rho_M(C)=0$, we have
$\rho^2(L)=(n_2^2-1)^2\rho^2(B)/(n^2-1)^2=0.0621\rho^2(B)=0.00002$
and $\rho_M(L)=(n_2^2-1)\rho_M(B)/(n^2-1)=0.2493\rho_M(B)=0.0191$.

A more general result than Proposition \ref{proposition4} can be obtained
if $A$ and $D$ are nearly orthogonal
and at least one of the two, $A^TC =0$ and $B^T D=0$,
approximately holds. However, besides being very complicated,
such a general result
does not greatly enhance our capability of
constructing nearly orthogonal Latin hypercubes
as the orthogonality
of $A$ and $D$ and that between $A$ and $C$
is much easier to achieve than the orthogonality of $B$ and $C$.
Our result as in Proposition \ref{proposition4} makes a more focused
presentation. Lin (\citeyear{L08}) also contains a table of small,
nearly orthogonal Latin hypercubes, based on which we can
construct large nearly orthogonal Latin hypercubes via Proposition \ref
{proposition4}.
\end{Example}

%s5 ###
\section{Concluding remarks}\label{sec5}

We have presented a general method of construction
for orthogonal, nearly orthogonal and cascading Latin hypercubes.
The method uses small designs to build
large designs. It turns out that some appealing properties in small
designs can be carried over to large designs. We have also obtained
a result on the existence of orthogonal Latin hypercubes.
The power of the general method is further
enhanced by the methods of stacking.
Although our methods are motivated by computer experiments,
they are potentially useful for constructing other designs such as
permutation arrays which are widely applied to data transmission
over power lines [see Colbourn,
Kl{\o}ve and Ling (\citeyear{CKL04}) and the reference therein].

Many researchers are increasingly interested in using polynomial
models for computer experiments though Gaussian process models
are still very popular. Polynomials are attractive because
they allow gradual building of a suitable model
by starting with simple linear terms and
then gradually introducing higher-order terms.
Orthogonal and nearly orthogonal Latin hypercubes are directly useful
when polynomial models are considered. If one insists on using
Gaussian-process models, orthogonality and near orthogonality
can be viewed as stepping stones to space-filling designs.
This is because a good space-filling
design must be orthogonal or nearly so as the design points
when projected on to two dimensions should be uniformly scattered.
Thus the search for space-filling designs can be restricted
to orthogonal and nearly orthogonal designs instead of all designs.
A rich class of orthogonal and nearly orthogonal Latin hypercubes
can be obtained by considering a generalization of the construction
method in this paper. The generalization makes use of an idea
in Bingham, Sitter and Tang (\citeyear{BST09}) [for more details,
we refer to Lin (\citeyear{L08})]. It is part of our research plan
to write a paper on this topic in the future.

\begin{appendix}\label{app}
\section*{Appendix}

\begin{pf*}{Proof of Lemma \protect\ref{lemma1}}
We provide a proof under (a) in condition (ii) of Lemma \ref{lemma1}.
The proof is essentially the same if condition (b) is met.
For design $L$ in (\ref{eq1}) to be a Latin hypercube, we need
to show that each column of $L$ is
a permutation of
$-(n-1)/2, -(n-3)/2, \ldots, (n-3)/2, (n-1)/2$ where $n=n_1n_2$.
Without loss
of generality, we will prove that this is the case for the first
column of design $L$. For ease in notation, let
$(a_1, \ldots, a_{n_1})^T$,
$(b_1, \ldots, b_{n_2})^T$,
$(c_1, \ldots, c_{n_1})^T$ and
$(d_1, \ldots, d_{n_2})^T$ be the first columns
of $A$, $B$, $C$ and $D$, respectively.
Then the entries of the first column of $L$ are given by
%
%e5.1 ###
\begin{equation}\label{eq:abcd}
a_i b_j + n_2 c_i d_j
\qquad\mbox{where $i=1,\ldots, n_1$ and $j=1,\ldots,n_2$.}
\end{equation}
As $C$ is a Latin hypercube, we have that $c_1, \ldots, c_{n_1}$
are a permutation of $-(n_1-1)/2, -(n_1-3)/2, \ldots, (n_1-3)/2, (n_1-1)/2$.
For any given odd $u$ such that $ 1 \leq u \leq n_1$, consider the two
distinct levels, $- (n_1 - u)/2$ and $ (n_1 - u)/2$, of~$C$.
(The two levels may be the same level $0$ when $n_1$ is odd. This
simple case
will be dealt with later.)
For this given $u$, let $i$ and $i'$ be the unique indices such that
$c_i = (n_1 - u)/2$ and $c_{i'} = - (n_1 - u)/2$.
As $d_j = \pm1$, the two numbers $c_i d_j$ and $c_{i'} d_j$
must always have opposite signs and thus always give
the two points $- (n_1 - u)/2$ and $(n_1 - u)/2$ on the real line.
Therefore, the two numbers $n_2 c_i d_j$ and $n_2 c_{i'} d_j$
always give the two points $- n_2 (n_1 - u)/2$ and $n_2 (n_1 - u)/2$
for any $j=1, \ldots, n_2$. By condition (a), we have that
$a_i = a_{i'}$. Since $B$ is a Latin hypercube of $n_2$ runs,
we have that $b_1, \ldots, b_{n_2}$ are a permutation
of $- (n_2 -1)/2, -(n_2-3)/2, \ldots, (n_2-3)/2, (n_2 -1)/2$.
As $a_i = \pm1$, we have that $a_i b_1, \ldots, a_i b_{n_2}$
are also a permutation of
$- (n_2 -1)/2, -(n_2-3)/2, \ldots, (n_2-3)/2, (n_2 -1)/2$.
Since $a_{i'} = a_i$,
this shows that the $2n_2$ points given by
$a_i b_j + n_2 c_i d_j$ and
$a_{i'} b_j + n_2 c_{i'} d_j$
for $j =1, \ldots, n_2$ can be divided into two sets of $n_2$ points
with the
first set of $n_2$ points given by
$- n_2 (n_1 - u)/2 + b_j$ for $j=1, \ldots, n_2$
and the second set of $n_2$ points given by
$ n_2 (n_1 - u)/2 + b_j$ for $j=1, \ldots, n_2$.
The $n_2$ points $- n_2 (n_1 - u)/2 + b_j$ for $j=1, \ldots, n_2$
are centered at $- n_2 (n_1 - u)/2 $, and equally spaced with
two adjacent points separated by an interval of length one.
A similar remark can be made about the other set of $n_2$ points.
We note that if $u=n_1$ when $n_1$ is odd, for the unique $i$ with $c_i =0$,
the $n_2$ numbers
$
a_i b_j + n_2 c_i d_j = a_i b_j
$
for $j=1, \ldots, n_2$
are simply the set of $b_j$s for $j=1, \ldots, n_2$.
By allowing the odd $u$ to vary in the range $ 1 \leq u \leq n_1$,
we see that the $n_1 n_2$ numbers in (\ref{eq:abcd}) are precisely
these $n$ points, $-(n-1)/2, -(n-3)/2, \ldots, (n-3)/2, (n-1)/2$, where
$n=n_1 n_2$.
The proof is complete.
\end{pf*}
\begin{pf*}{Proof of Theorem \protect\ref{theorem2}}
The sufficiency part of Theorem \ref{theorem2} can be proved directly
which involves the construction of an orthogonal Latin hypercube
of $n$ runs with $m \geq2$ factors for any $n$ that does not
have form $4k+2$. We omit this part of the proof as
the existence result also follows from Proposition \ref{proposition3}
in Section~\ref{sec33} when we establish a lower bound on the maximum
number of factors in an orthogonal Latin hypercube.

It remains to show that
there does not exist an orthogonal Latin hypercube of
$n=4k+2$ runs with $m \geq2$ factors.
Now suppose that such an orthogonal Latin hypercube exists,
and let $a=(a_1, \ldots, a_n)^T$ and $b=(b_1,\ldots, b_n)^T$ be its
two columns. Then we have that
both $a$ and $b$ are permutations of $\{1/2,3/2, \ldots,
(n-1)/2,-1/2,-3/2,\ldots,-(n-1)/2\}$.
Note that $\sum_{i=1}^n a_i =
0$, $\sum_{i=1}^n b_i=0$. Without loss of generality, we assume that
$a = (1/2,3/2$,\break $\ldots,
(n-1)/2,-1/2,-3/2,\ldots,-(n-1)/2)^T$. In other words,
we have $a_i=-a_{i+n/2}=(2i-1)/2$.
Since $a$ and $b$ are orthogonal, we have that
$\sum_{i=1}^n a_i b_i= 2^{-1}\sum_{i=1}^{n/2}[(2b_i)i-(2b_{i+n/2})(i-1)]=0$.
Note that
both $2b_i$ and $2b_{i+n/2}$ are odd, $i=1, \ldots, n/2$. The
quantity $(2b_i)i-(2b_{i+n/2})(i-1)$ must be odd as $(2b_i)i$ and
$(2b_{i+n/2})(i-1)$ cannot be both even or both odd. In addition,
$n/2$ must be odd. It is obvious that the addition or subtraction
among an odd number of odd integers gives an odd integer. This
leads to a contradiction.
\end{pf*}
\begin{pf*}{Proof of Proposition \protect\ref{proposition4}}
Parts (i) and (ii) can be obtained by noting that
\begin{eqnarray*}
L^TL&=& (A \otimes B + \gamma C \otimes D)^T(A \otimes B
+ \gamma C \otimes D)\\
&=& (A^TA)\otimes(B^TB)+\gamma(A^TC)\otimes(B^TD)\\
&&{}
+ \gamma(C^TA)\otimes(D^TB) +\gamma^2 (C^TC)\otimes(D^TD)\\
&=& n_1I_{m_1}\otimes(B^TB)+n_2^2(C^TC)\otimes(n_2I_{m_2}),
\end{eqnarray*}
where $I_{m_1}$ and $I_{m_2}$ are identity matrices of
size $m_1$ and $m_2$, respectively. The second step follows
by the properties of the Kronecker product given in (\ref{eq:kron}).
The last step is due to the orthogonality of $A$ and $D$,
either of the conditions $A^TC=0$ and $B^TD=0$, and $\gamma=n_2$.
In addition, for an $n \times m$ Latin hypercube $L$, the sum of squares
of the elements in each of its columns is $n(n^2-1)/12$.
Thus the $m \times m$ correlation matrix among the $m$ columns of $L$
is given by $[n(n^2-1)/12]^{-1}L^TL$. Based on the elements in the
correlation matrix, $\rho^2(L)$ and $\rho_M(L)$ can be computed
in the following way:
\begin{eqnarray*}
\rho^2(L)&=&\bigl(m_1n_1^2m_2(m_2-1)[n_2(n_2^2-1)/12]^2\rho^2(B)\\
&&\hspace*{2pt}{}
+n_2^6m_2m_1(m_1-1)[n_1(n_1^2-1)/12]^2\rho^2(C)\bigr)\\
&&\hspace*{0pt}{}\times\bigl(m_1m_2(m_1m_2-1)
[n(n^2-1)/12]^2\bigr)^{-1}\\
&=&\frac{(m_2-1)(n_2^2-1)^2\rho^2(B)+n_2^4(m_1-1)(n_1^2-1)^2\rho^2(C)}
{(m_1m_2-1)(n^2-1)^2}
\end{eqnarray*}
and $\rho_m(L)$ is the larger value between $n_1n_2[(n_2^2-1)/12]\rho
_M(B)/[n(n^2-1)/12]$ and
$n_2^3n_1[(n_1^2-1)/12]\rho_M(C)/[n(n^2-1)/12]$. With the definition of
$w_1$, $w_2$, $w_3$ and $w_4$, we complete the proof.
\end{pf*}
\end{appendix}

\section*{Acknowledgments}
The authors thank the Associate Editor and referees for
their helpful comments.

%
%accubitu suo, nardus mea dedit odorem suavitatis. Quoniam confortavit
%seras portarum tuarum, benedixit filiis tuis in te. Qui posuit fines
%tuos}

%
\printaddresses

\end{document}